\documentclass[12pt]{article}
\usepackage[utf8]{inputenc}
\usepackage{amssymb,amsthm,bm,graphicx,xcolor, tipa,tikz}
\usepackage{mathtools}
\usepackage[T1]{fontenc}

\newtheorem{thm}{Theorem}[section]

\newtheorem{conj}[thm]{Conjecture}

\newtheorem{rem}[thm]{Remark}
\newtheorem*{cor*}{Corollary}
\newtheorem{lem}[thm]{Lemma}
\newtheorem{fact}[thm]{Fact}

\newtheorem{cla}[thm]{Claim}

\newtheorem{defi}[thm]{Definition}

\theoremstyle{definition}

\theoremstyle{remark}

\newtheorem*{rek*}{Remark}

\newtheorem*{cla*}{Claim}

\newcommand{\Z}{\mathbb{Z}}
\newcommand{\C}{\mathcal{C}}
\newcommand{\M}{\mathcal{M}}

\newcommand{\R}{\mathbb{R}}

\newcommand{\bp}{\begin{proof}}
\newcommand{\ep}{\end{proof}}
\newcommand{\eq}[1]{\begin{align*}#1\end{align*}}
\newcommand{\eqn}[1]{\begin{align}#1\end{align}}

\newcommand{\fracflat}[2]{#1/#2}

\DeclarePairedDelimiter\ceil{\lceil}{\rceil}

\title{Coprime Mappings and Lonely Runners}
\author{Tom Bohman\thanks{This work was supported by a grant from the Simons
Foundation (587088, TB)} \, and Fei Peng}
\date{\today}

\begin{document}

\maketitle

\begin{abstract}
    For $x$ real, let $ \{ x \}$ be the fractional part of $x$ (i.e. $\{x\} = x - \lfloor x \rfloor $).  
    The lonely runner conjecture can be stated as follows: for any $n$ positive integers $ v_1 < v_2 < \dots < v_n $ 
    there exists a real number $t$ such that $ 1/(n+1) \le \{ v_i t\} \le n/(n+1) $ for $ i = 1, \dots, n$.  In this 
    paper we prove that if $ \epsilon >0 $ and $n$ is sufficiently large (relative to $\epsilon$) then
    such a $t$ exists for any collection of positive integers $ v_1 < v_2 < \dots < v_n$ such that 
    $ v_n < (2-\epsilon)n$.  This is an approximate version of a natural next step for the study of the lonely 
    runner conjecture suggested by Tao.  
    
    The key ingredient in our proof is a result on coprime mappings.  Let $A$ and $B$ be sets of integers. 
    A bijection $ f:A \to B$ is a coprime mapping if $ a $ and $f(a)$ are coprime for every $ a \in A$.
    We show that if $A,B \subset [n]$ are intervals of length $2m$ where $ m = e^{ \Omega({(\log\log n)}^2)}$ then there
    exists a coprime mapping from $A$ to $B$.  We do not believe that this result is sharp.

\end{abstract}

\section{Introduction}
Suppose $n$ runners are running on a circular track of circumference 1.  It is not a race.  The runners all start at the same point on the track and at the same time, and each runs at their own distinct constant speed.  We say that a runner is lonely at time $t$ if the distance (along the track) to the nearest of the other runners is at least $1/n$.  The lonely runner conjecture asserts that every runner is lonely at some point in time.  This problem originally arose in the context of diophantine approximations and view obstruction problems \cite{wills1967zwei} \cite{cusick1973view}. (The poetic formulation given here is due to Goddyn \cite{bienia1998flows}.)  It is easier to work with the following restatement of the conjecture, which we obtain by subtracting the speed of one runner from all speeds (then one of the runners is `standing still').  In the original problem speeds are real valued, but it is known that the general problem can be reduced to case where all speeds are integers (see \cite{bohman2001six}).  So we henceforth consider only integer speeds.
For $x$ real, let $ \{ x \}$ be the fractional part of $x$ (i.e. $\{x\} = x - \lfloor x \rfloor $). 


\begin{conj}[Lonely Runner Conjecture] \label{runnerdream}
 For any $n$ positive integers $ v_1 < v_2 < \dots < v_n $, \eqn{\label{constraint}\exists t\in \R\ \text{ such that }\ 1/(n+1) \le \{ v_i t\} \le n/(n+1)  \ \text{ for } \  i = 1, \dots n.}
\end{conj}

\noindent
There are examples of sets of speed which ``almost'' break Condition \eqref{constraint}. 

\begin{defi}
\label{tinsdef}Positive integers $ v_1 < v_2 < \dots < v_n $ are said to be a {\bf tight instance} for the lonely runner conjecture if Condition \eqref{constraint} holds, but only with equality.  In other words, the instance $ v_1 < v_2 < \dots < v_n $ is tight if (\ref{constraint}) holds 
and  there does not exist $ t\in \R$ such that $ 1/(n+1) < \{ v_i t\} < n/(n+1)$  for $i = 1, \dots n$.
An instance that is neither a counterexample nor tight is a {\bf loose instance} of the lonely runner conjecture. An instance is loose if
\eqn{\label{loose}\exists t\in \R\ \text{ such that }\ 1/(n+1) < \{ v_i t\} < n/(n+1)  \ \text{ for } \  i = 1, \dots n.}
\end{defi}

\noindent
The canonical example of a tight instance is $(1,2,\dots,n)$.
Tight instances were studied by Goddyn and Wong \cite{Goddyn2006TIGHTIO}. They showed that the canonical tight instance can be modified to create another tight instance by accelerating a speed that is slightly less than $n$ -- and satisfies certain number theoretic conditions -- by a suitable integer factor.  For example,  $(1,2,3,4,5,7,12)$  is a tight instance.  They also showed that small sets of speeds in the canonical instance that satisfy these conditions can be simultaneously accelerated to produce tight instances.

Tao \cite[Proposition 1.5]{tao2017remarks} showed that \eqref{constraint} holds if $v_1,\dots,v_n\le 1.2n$. He also suggested that proving the conjecture holds 
for $ v_1, \dots, v_n \le 2n$ is a natural target: the condition $v_1,\dots,v_n\le 2n$, unlike $1.2n$, allows multiple tight instances. So the desired statement would prove the conjecture for instances that are in the vicinity of tight instances. We prove an approximate version of this target:
\begin{thm}
\label{thm:run}
There exists a constant $\C$ such that for sufficiently large $n$, if $n<v_n\le 2n-\exp(\C\cdot(\log\log n)^2)$, then positive integers $v_1<v_2<\dots<v_n$ are a loose instance for the lonely runner conjecture.
\end{thm}
\noindent
Unfortunately, as seen in the theorem statement, the underlying objective to separate tight instances from counterexamples is not achieved.

The key ingredient in our proof of Theorem~\ref{thm:run} is inspired by coprime mappings.
\begin{defi}
If $ A, B$ are sets of integers then a bijection $ f: A \to B$ is a coprime mapping if $ a$ and $ f(a)$ are coprime for every $a \in A$.
\end{defi}
\noindent
Initial interest in coprime mappings was focused on the case $A = [n]$. 
D.J. Newman conjectured that for all $n\in\Z^+$ there is a coprime mapping between $[n]$ and any set of $n$ consecutive integers.
This conjecture was proved by Pomerance and Selfridge \cite{pomerance1980proof} (after Daykin and Baines \cite{daykin_baines_1963} 
and Chvátal \cite{chvatal} established special cases of the conjecture.) Robertson and Small \cite{robertsonsmall} determined when 
a coprime mapping exists between $A=[n]$ (or $A=\{1,3,5,\dots,2n-1\}$) and an $n$-term arithmetic progression (AP). 

More recently there has been interest in coprime mappings where neither $A$ nor $B$ contains 1.  Note that if $A = \{2, \dots, n+1 \}$, the integer $s$ is the product of all primes that are at most $n+1$,  
and $ s \in B $ then there is no coprime map from $A$ to $B$.  So we must place some restriction on the set $B$ if we consider sets 
$A$ that do not contain 1.
Larsen et al \cite{Larsen2017CoprimeMO} considered sets of adjacent intervals of integers. They conjectured that
if $ 1 \le \ell < k$ and $ k \neq 3$ then there is coprime mapping from $A=\{\ell+1,\dots,\ell+k\}$ to $ B=\{\ell+k+1,\dots,\ell+2k\}$. 

The application of coprime mappings that we use in the context of the lonely runner conjecture requires the further generalization to the case
that $A$ and $B$ are not adjacent.  For each positive integer $n$ we define the number $ f(n)$ to be the smallest integer
such that for all $2m \ge f(n)$ there is a coprime mapping between every pair of intervals $A,B \subset [n]$ with $|A|=|B|=2m$.\footnote{We force the cardinality to be even, because when the cardinality is odd, $A$ and $B$ can both have a majority of even numbers, making a coprime mapping impossible.}
Note that the example we give above
establishes the bound
\eq{ f(n) > (1-o(1))\log n. } 
Note further that the Conjecture of Larsen et al requires only a linear upper bound on $f(n)$. We establish a stronger asymptotic bound.
\begin{thm}
\label{thm:co}
$f(n) = \exp(O((\log\log n)^2)). $
\end{thm}
\noindent
We do not believe that this result is sharp. Indeed, we conjecture that $ f(n)$ is at most polylogarithmic in $n$. We also resolve the conjecture of Larsen et al.
\begin{thm}
\label{tnconj}
If $ 0 \le \ell < k$ and $ k \ge 4$ then there is coprime mapping from $A=\{\ell+1,\dots,\ell+k\}$ to $\ B=\{\ell+k+1,\dots,\ell+2k\}$. 
\end{thm}
\noindent
Note that there is no coprime mapping from $A = \{2,3,4\} $ to $ B = \{ 5,6,7\}$.  Thus, some condition on $ k$ is required.

The remainder of this paper is organized as follows. In the next Section we prove our central result, which can be viewed as a number theoretic version of Hall's condition.  Theorem~\ref{thm:co} follows immediately from the central result.  In Section~3 we use the central result to prove Theorem~\ref{thm:run}. In the final section we prove the conjecture of Larsen et al regarding coprime mappings between adjacent intervals; that is, we prove Theorem~\ref{tnconj}.

\section{The central result}

Our central result is as follows:
\begin{thm}\label{evencase}
There exists a constant $\C$ such that the following is true for sufficiently large $n$. If $I,J\subseteq [n]$ are both sets of $2m$ consecutive integers and $2m\ge\exp(\C\cdot(\log\log n)^2)$, then for all subsets $S\subseteq I$, $T\subseteq J$ that satisfy $|S|+|T|\ge 2m$, exactly one of the following happens:
\begin{enumerate}
    \item $S=\emptyset$;
    \item $T=\emptyset$;
    \item $S=I\cap 2\Z$ and $T=J\cap 2\Z$;
    \item there exist $s\in S$, $t\in T$ that are coprime.
\end{enumerate}
In particular, if $|S|+|T|>2m$ then there is a coprime pair.
\end{thm}

\begin{rem}
In order to apply Hall's Theorem to establish the existence of a coprime mapping between $I$ and $J$, it suffices to show that every pair of sets $ S \subseteq I, T \subseteq J$ such that $ |S| + |T| \ge 2m+1$ contains a coprime pair $s,t$ such that $s \in S$ and $t \in T$. Thus, Theorem~\ref{thm:co} follows immediately from Theorem~\ref{evencase}. Note that Theorem~\ref{evencase} is stronger than necessary for this purpose as it treats the case $ |S| + |T| = 2m$. This case is needed for the application to the lonely runner problem. (I.e. for the proof of Theorem~\ref{thm:run}.)
\end{rem}

\noindent
We now turn to the proof of Theorem~\ref{evencase}.
Lemma~\ref{interlem} is the core of the proof.  It has a weaker condition (APs) and a weaker result (2-coprime) than Theorem~\ref{evencase}. For the remainder of this section we assume $ I,J \subset [n]$ are APs of cardinality $m$ with common difference 1 or 2.
\begin{defi}
Two numbers $s,t\in\Z^+$ are said to be 2-coprime if no prime other than 2 divides them both. E.g., (3,4), (12,16).
\end{defi}

\begin{lem}\label{interlem}
There exists a constant $\C$ such that the following is true for sufficiently large $n$. 
If $m \ge \exp(\C\cdot(\log\log n)^2)$ then for all
nonempty subsets $S\subseteq I$, $T\subseteq J$ such that $|S|+|T|\ge m$ there exist $s\in S,\ t\in T$ that are 2-coprime.
\end{lem}

\noindent
The main ingredients of the proof of Lemma~\ref{interlem} are the following two lemmas and one fact.

\begin{lem}\label{rgcoprime}
Let $S\subseteq I$, $T\subseteq J$ be nonempty subsets such that $|S|+|T|\ge m$, and let $r=m/|S|$. If $r\ge16$, $m>5\log(n)^{\log_2(2r)}$, and $n$ is sufficiently large 
then there exist $s\in S,\ t\in T$ that are 2-coprime.
\end{lem}

\begin{lem}\label{rlcoprime}
Let $S\subseteq I$, $T\subseteq J$ be nonempty subsets such that $|S|+|T|\ge m$, and let $r=m/|S|$. If $ 2\le r\le16$, $m>\log(n)^3$ and $n$ is sufficiently large then there exist $s\in S,\ t\in T$ that are 2-coprime.
\end{lem}

\begin{fact}[See \cite{10.2307/24489279}]\label{coprimeconsec}
There are positive constants $c_1,c_2$ such that for all $n\in \Z^+$, among any sequence of $c_1\cdot \omega(n)^{c_2}$ consecutive integers, there is at least one that is coprime to $n$. (Here $\omega(n)$ is the number of distinct prime divisors of $n$.)
\end{fact}
\noindent
We now prove Lemma \ref{interlem}, assuming Lemma~\ref{rgcoprime} -- Fact~\ref{coprimeconsec}. 
We will prove Lemmas \ref{rgcoprime} and \ref{rlcoprime} immediately after. We end the section with the proof of Theorem~\ref{evencase}.

\bp[Proof of Lemma \ref{interlem}]

Set $\C =2c_2\log_2(e)$, where $c_2$ is the constant in Fact~\ref{coprimeconsec}.
Without loss of generality, assume that~$|S|+|T|=m$. Since the roles of $S$ and $T$ are the same, one can assume that $|S|\le m/2$. That is, $r:=m/|S|\ge 2$. Note that Lemma~\ref{interlem} follows immediately from either Lemma~\ref{rgcoprime} or Lemma~\ref{rlcoprime} (depending on the value of $r$) unless 
%
\eqn{ \label{larger} \exp(\C\cdot(\log\log n)^2) \le m < 5\log(n)^{\log_2(2r)} .}
(Note that we clearly have $ \exp(\C\cdot(\log\log n)^2) > \log(n)^3$ for $n$ sufficiently large.) It remains to prove Lemma~\ref{interlem} when $r$ and $m$ satisfy (\ref{larger}). 
To this end, we first observe that if we replace coprime with 2-coprime, then we can extend Fact~\ref{coprimeconsec} to APs with common difference 2.

\begin{cla*}
 For all $n\in \Z^+$, any integer AP with common difference 1 or 2 and at least $c_1\cdot \omega(n)^{c_2}$ terms contains at least one term that is 2-coprime to $n$. (Using the same constants as in Fact~\ref{coprimeconsec}.)
\end{cla*}
\begin{quote}
\bp[Proof of Claim] 
Assume without loss of generality that $n$ is odd.
Say the AP is $a_1,\dots,a_\ell$. If the AP contains only even numbers, consider $a_1/2,a_2/2,\dots,a_\ell/2$; if only odd numbers, consider $(a_1+n)/2,(a_2+n)/2,\dots,(a_\ell+n)/2$. In any case, we have a sequence of $\ell$ consecutive integers, where the $i$-th number is coprime to $n$ if and only if $a_i$ is coprime to $n$. Fact~\ref{coprimeconsec} implies that the new sequence has a number coprime to $n$.
\ep
\end{quote}

Now consider $s\in S$. By the Claim, among any $c_1\cdot \omega(s)^{c_2}$ consecutive terms of $J$, at least one is 2-coprime to $s$. Thus, $J$ has at least $\lfloor m/(c_1\omega(s)^{c_2})\rfloor$ terms that are 2-coprime to $s$. If any of these terms are in $T$ then we have the desired 2-coprime pair, so we may assume for the sake of contradiction that they
all in $J\setminus T$. 
Note that, $\omega(s)\le \log_2(s)\le \log_2(n)$, and, for sufficiently large $n$, we have \eq{\frac{m}{c_1\omega(s)^{c_2}}\ge\frac{\exp(\C\cdot(\log\log n)^2)}{c_1\log_2(n)^{c_2}}.} 
As this quantity is arbitrarily large for large $n$, the floor function has negligible effect, and it follows that we have 
\eq{ r = \frac{m}{|S|} = \frac{m}{|J \setminus T|} \le \frac{m}{\left\lfloor\frac{m}{c_1\omega(s)^{c_2}}\right\rfloor} \le (1+o(1)) c_1\log_2(n)^{c_2} .}
Then, again appealing to (\ref{larger}), we have \eq{ m < 5\log(n)^{\log_2(2r)} \le 5\log(n)^{(c_2+o(1))\log_2\log_2n}<\exp(\C\cdot(\log\log n)^2)\le m.}
This is a contradiction.
\ep

We now prove Lemmas \ref{rgcoprime} and \ref{rlcoprime}. The key idea is to count the non-coprime pairs in $S\times T$ by summing up $|S\cap p\Z||T\cap p\Z|$ over primes $p$ greater than 2. Note that if $S$ and $T$ are random subsets of $I$ and $J$, respectively, than we expect to have 
\eq{\sum_{p >2} |S\cap p\Z||T\cap p\Z| \approx \sum_{p>2} \frac{|S|}{p} \frac{|T|}{p} = |S \times T| \sum_{p >2} \frac{1}{p^2} \approx 0.2 |S \times T|,}
and there are many 2-coprime pairs.
Of course, we have to complete the proof for all sets $S$ and $T$ (rather than just random ones) and we do this by `zooming in' on primes $p$ for which $ |S \cap p\Z|$ is large.  We `zoom in' by looking at $S\cap p\Z$ and $T\setminus p\Z$ instead of $S$ and $T$, and we iterate this process if necessary.

Before the proof, we state a simple approximation that we use throughout this Section.
\begin{lem}
\label{lem:simpleapprox}
Let $A$ be an integer AP with common difference $d$, and suppose $P\in\Z^+$ is coprime with $d$. Then for all $R\subseteq [P]$, \eqn{\label{absrem}\frac{|A\cap (R+P\Z)|}{|R|}-\frac{|A|}{P}\in(-1,1).}
As a result, if $|A|/P\ge \delta^{-1}$ for some $\delta>0$, then \eqn{\label{relrem}\frac{|A\cap (R+P\Z)|}{|A||R|/P}\in(1-\delta,1+\delta).}
\end{lem}
\bp
Among every $P$ consecutive terms of $A$ (a ``chunk''), exactly $|R|$ of them are in $A\cap (R+P\Z)$. $\lfloor |A|/P\rfloor$ disjoint chunks can cover a subset of $A$ and $\lceil |A|/P\rceil$
chunks can cover a superset of $A$, so $|A\cap (R+P\Z)|$ is between $\lfloor |A|/P\rfloor |R|$ and $\lceil |A|/P \rceil |R|$.
\ep

\bp[Proof of Lemma \ref{rgcoprime}]
Let $S_0=S,\ T_0=T,\ I_0=I,\ J_0=J$. Let \eqn{\label{capmdefrg}M={\left(\frac{m}{5}\right)}^{1/\log_2(2r)}.} Note that by assumption, \eq{M>\left(\log(n)^{\log_2(2r)}\right)^{1/\log_2(2r)}=\log(n).}
For each $i\in\Z^+$, if there exists prime $p_i\notin\{2,p_1,\dots,p_{i-1}\}, p_i\le M$ such that \eq{
\frac{\left|S_{i-1}\cap p_i\Z\right|}{\left|I_{i-1}\cap p_i\Z\right|}\ge 2\frac{|S_{i-1}|}{|I_{i-1}|},
} define \eq{S_i=S_{i-1}\cap p_i\Z,\  T_i=T_{i-1}\setminus p_i\Z,\ I_i=I_{i-1}\cap p_i\Z,\  J_i=J_{i-1}\setminus p_i\Z.}
Let's say $k$ is the last index where these are defined. For every $0\le i\le k$ define $P_i=\mathtt{Primes}\setminus\{2,p_1,\dots,p_i\}$ for convenience. Now we have \eqn{
\label{skconc} &\frac{|S_i|}{|I_i|} \ge 2\frac{|S_{i-1}|}{|I_{i-1}|} \hskip1cm & \text{ for } i =1, 
\dots k, \text{ and }\\
\label{gkconc} & \frac{\left|S_k\cap p\Z\right|}{\left|I_k\cap p\Z\right|} < 2\frac{|S_k|}{|I_k|} \hskip3mm  &
\forall p\in P_k \text{ such that } p\le M.
}
Here, \eqref{skconc} implies that \eq{\frac{|S_k|}{|I_k|}\ge 2^k \frac{|S_0|}{|I_0|}=\frac{2^k}r.}
Since $S_i$ is a subset of $I_i$, we have $ |S_k|/ |I_k| \le 1$, and hence $k\le \log_2(r)$.
Define $\Gamma=\prod_{i=1}^k p_i$. Note that \eqn{\label{mgamma}M\Gamma=M\prod p_i\le M\cdot M^k \le M^{\log_2(2r)}=\frac m{5}.}

We establish some estimates for $ |I_k|, |J_k| $ and $ |I_k \cap p \Z|, |J_k \cap p\Z|$ for $ p \in P_k$ such that $p < M$.
Since $I_k=I\cap \Gamma\Z$ and $|I|/\Gamma\ge 5M$, by \eqref{relrem} we have \eqn{\label{ikest}|I_k|\in ( 1 \pm 1/M)m/\Gamma.}
$I_k$ is once again an AP. For $p\in P_k$ such that $p\le M$, $p$ does not divide the common difference of $ I_k$ and $|I_k|/p\ge |I_k|/M>(1-1/M)m/M\Gamma\ge 5-5/M$. Hence by \eqref{relrem} \eqn{\label{ikpest}|I_k\cap p\Z|\in \left(\left(1-\frac1{5-5/M}\right)\frac{|I_k|}{p},\  \left(1+\frac1{5-5/M}\right)\frac{|I_k|}{p}\right).}
Similarly, $J_k$ and $J_k\cap p\Z$ can be bounded from both sides. Since $J_k=J\setminus p_1\Z\setminus \dots\setminus p_k\Z$, and $|J|/\Gamma\ge 5M$, by \eqref{relrem} we have
\eq{|J_k|\in ( 1 \pm 1/M)m\cdot \prod(p_i-1)/\Gamma.}
Define $\Phi=m\prod(p_i-1)/\Gamma$. Then \eqn{\label{jkest}|J_k|\in( 1\pm 1/M) \Phi.}
For $J_k\cap p\Z$, we apply \eqref{relrem} with $A=J$, $P=\Gamma p$, $R=([\Gamma p]\cap p\Z)\setminus p_1\Z\setminus\dots\setminus p_k\Z$. Since $|J|/\Gamma p\ge 5$, we have \eqn{\label{jkpest}\forall p\in P_k,p\le M,\ |J_k\cap p\Z|&\in \left(0.8m|R|/|P|,\ 1.2m|R|/|P|\right)\nonumber\\&=(0.8\Phi/p, 1.2\Phi/p).}

Now consider $\Phi$. Note that 
\eq{\Phi =m\cdot \prod_{i=1}^k\left(1-\frac1{p_i}\right)
& \ge m\cdot \prod_{i=1}^{k}\left(1-\frac1{i+2}\right)
=2m/(k+2) \\
&  \ge 2m/\log_2(4r)
> \frac{2r |S| }{\log_2(4r)}
> \frac{16}{3}|S| \ \ \ \ \ (\text{as } r\ge 16).}
That is, $|S|<\fracflat{3\Phi}{16}.$
Since $J_k\setminus T_k\subseteq J\setminus T$, whose cardinality is at most $|S|$, by \eqref{jkest}, we have
\eqn{\label{tkjkest}\frac{|T_k|}{|J_k|}=1-\frac{|J_k\setminus T_k|}{|J_k|}&>1-\frac{3\Phi/16}{(1-1/M)\Phi}=\frac{13/16-1/M}{1-1/M}>\frac34;
\\ \label{tkest}|T_k|&> \frac{13/16-1/M}{1-1/M}|J_k|\ge (13/16-1/M)\Phi\ge\frac{3\Phi}4
.}

With these estimates in hand, we consider the quantity \eq{\lambda(p)=\frac{\left|S_k\cap p\Z\right|}{|S_k|}\frac{\left|T_k\cap p\Z\right|}{|T_k|}}
for every odd prime $p$. Note that the set of pairs $(s,t) \in S_k \times T_k$ that are {\bf not} 2-coprime because $p$ divides both $s$ and $t$ is $ (S_k\cap p\Z) \times (T_k\cap p\Z) $. Therefore, it follows from pigeonhole that if $\sum_{p\in P_0}\lambda(p)<1$ then there is a pair $(s,t)$ that is none of these sets. Such a pair is 2-coprime.  Thus, it suffices to show $\sum_{p\in P_0}\lambda(p)<1$. In order to estimate $ \lambda(p)$ we divide into cases:

\begin{itemize}
    \item \textbf{When $p\in \{p_1,\dots,p_k\}$:} By definition, $T_k$ has already excluded multiples of such $p$, so $\lambda(p)=0$.
    \item \textbf{When $p\in P_k,p\le M$:} By \eqref{gkconc} and \eqref{ikpest}, \eqn{\label{skpskest}\frac{\left|S_k\cap p\Z\right|}{|S_k|}< 2\frac{\left|I_k\cap p\Z\right|}{|I_k|}< \left(1+\frac1{5-5/M}\right)\frac{2}{p}<\frac{5}{2p}.} By \eqref{jkpest} and \eqref{tkest}, \eq{\frac{\left|T_k\cap p\Z\right|}{|T_k|}\le\frac{\left|J_k\cap p\Z\right|}{|T_k|}<\frac{1.2\Phi/p}{3\Phi/4}=\frac{8}{5p}.}
    Hence, $\lambda(p)<\frac{4}{p^2}$.
    \item \textbf{When $p\in P_k, p>M$:} In this case, we don't have equally good individual bounds. Here we use the following simple observation: \eqn{\label{sskp}\sum_{p\in P_k, p>M} \left|S_k\cap p\Z\right| <\log_M(n) |S_k|.} (This follows from the fact that every number in $S_k$ is at most $n$, and so is divisible by less than $\log_M(n)$ primes greater than $M$.) 
    Also by \eqref{absrem} and \eqref{tkjkest}, \eq{\frac{|T_k\cap p\Z|}{|T_k|}\le\frac{|J_k\cap p\Z|}{|T_k|}<\frac{|J_k|/p+1}{|T_k|}<\frac{4}{3p}+\frac1{|T_k|}.} By \eqref{tkest} and \eqref{mgamma}, \eq{|T_k|>\frac{3\Phi}{4}=\frac{3m\prod(p_i-1)}{4\Gamma}\ge \frac{3m}{4\Gamma} > \frac{15M}4.}
    Since $p>M$, \eq{\frac{|T_k\cap p\Z|}{|T_k|}<\frac{4}{3M}+\frac{4}{15M}=\frac{8}{5M}.} In view of \eqref{sskp}, 
    \eq{\sum_{p\in P_k, p>M}\lambda(p)<\log_M(n) \cdot \frac{8}{5M} = \frac{ 8 \log(n)}{ 5 M \log M} .}
\end{itemize}

Summing up all cases, \eq{\sum_{p\in P_0}\lambda(p)&<\sum_{p\in P_k,p<M}\frac{4}{p^2}+\sum_{p\in P_k,p>M}\lambda(p)\nonumber\\ &< 4\left(P(2)-\frac14\right)+\frac{ 8 \log(n)}{ 5 M \log M} \nonumber\\ &< 4\left(P(2)-\frac14\right)+\frac{ 8 }{5 \log\log(n)}\tag{as $M>\log(n)$} \nonumber\\ & < 1,}
for $n$ sufficiently large,
where $P(2)=\sum_{p\text{ prime}}p^{-2}$ is the prime zeta function. (Note that $ \log \log (n) \ge 9 $ suffices.)
As this sum is less than 1, we conclude that some $s\in S_k\subseteq S$ and $t\in T_k\subseteq T$ are 2-coprime.\ep 

\vskip5mm

\bp[Proof of Lemma \ref{rlcoprime}]
Let $\alpha := 1/r = |S|/m\in [1/16, 1/2]$, $P_0=\mathtt{Primes}\setminus\{2\}$ and 
\eqn{\label{capmdefrl}M = m^{1/3} >\log(n).} 
Let $P_M=P_0\cap [M]$. For each $p\in P_0$, define $\alpha_p=|S\cap p\Z|/m\le \alpha$. 

We begin with some estimates.
Note that we have \eqn{\label{test}|T| \ge m-|S|=(1-\alpha)m.} For every $ p,q\in P_M$ $(p\neq q)$, $m/pq\ge m/M^2 = M$. By \eqref{relrem} we have,  
\eqn{\label{ipqest} |I\cap pq\Z| = ( 1 \pm  1/M) \frac{m}{pq} \hskip1cm  & \text {for all } p\neq q\in P_M ,\\
|(J\setminus p\Z)\cap q\Z| =  (1 \pm 1/M) \frac{(p-1)m}{pq}  \hskip1cm & \text{for all } p\neq q\in P_M, \text{ and }  \\
\label{ijpest}{|I\cap p\Z|}, {|J\cap p\Z|} \in (1 \pm 1/M) \frac{m}{p} \hskip1cm & \text{ for all } p\in P_M.
}
A consequence of \eqref{ijpest} and \eqref{test} is \eqn{\label{tpest}\forall p\in P_M,\ |T\setminus p\Z|\ge |T|-|J\cap p\Z|\ge (1-\alpha- ( 1 + 1/M)/p)m.}
As in Lemma \ref{rgcoprime} we have
\eqn{\label{ssps} \sum_{p\in P_0\setminus P_M} |S\cap p\Z| < |S| \log_M(n) &
;\\\label{sspqsp} \sum_{q\in P_0\setminus P_M} |(S\cap p\Z)\cap q\Z| <  | S \cap  p\Z| \log_M(n) &  \text{ for all } p\in P_0 \text{ such that } \alpha_p > 0.} 
For $p\in P_M$ such that $\alpha_p>0$, by \eqref{absrem} and \eqref{tpest}, we have \eqn{\label{tpqtp}\forall q\in P_0\setminus P_M,\ \frac{|(T\setminus p\Z)\cap q\Z|}{|T\setminus p\Z|} &\le \frac{|J\cap q\Z|}{(1-\alpha- (1 + 1/M)/p)m}\nonumber\\&\le\frac{m/q+1}{(1/6- 1 /(Mp))m}\nonumber\\&=\frac6M\frac{M/q+1/M^2}{1-6/Mp}\nonumber\\(\text{assuming }M\ge15,)\ \ &< \frac{7}{M}.}
(Remark: Note that it is crucial that we excluded 2 from the primes.) With these estimates in hand, we consider two cases: \\[0mm]

\noindent
{\bf Case (a)}:
There exists some $p\in P_M$ such that \eqn{\label{alphapboundn}\alpha_p> \frac{0.24}{p(1-3\alpha/2)} .}

\noindent    
Fix such $p\in P_M$, and for $q\in P_0$ define \eq{\lambda_1(q)=\frac{|(S\cap p\Z)\cap q\Z|}{|S\cap p\Z|}\frac{|(T\setminus p\Z)\cap q\Z|}{|T\setminus p\Z|}.} In view of \eqref{tpest}, $T\setminus p\Z$ is nonempty, and so we will have the desired 2-coprime pair if $\sum_{q\in P_0}\lambda_1(q)<1$. Note that when $q=p$, $\lambda_1(q)=0$. Thus, by \eqref{tpqtp}, 
\eqref{sspqsp}, 
\eqref{ijpest}, 
and \eqref{tpest}, we have  \eq{\sum_{q\in P_0}\lambda_1(q) &= \sum_{q\in P_M\setminus\{p\}}\lambda_1(q)+\sum_{q\in P_0\setminus P_M}\lambda_1(q)\\
&\le \sum_{q\in P_M\setminus\{p\}}\frac{|I\cap pq\Z|}{|S\cap p\Z|}\frac{|(J\setminus p\Z)\cap q\Z|}{|T\setminus p\Z|}+\left(\sum_{q\in P_0\setminus P_M}\frac{|(S\cap p\Z)\cap q\Z|}{|S\cap p\Z|}\right)\frac{7}{M}\\
&\le \sum_{q\in P_M\setminus\{p\}} \frac{ (1+1/M)m/pq}{\alpha_pm}\frac{(1+1/M)(p-1)m/pq}{(1 -1/M)(1-\alpha-1/p)m}\ +\  \log_M(n) \cdot \frac{7}{M}\\
&=\frac{ (1 + 1/M)^2}{ 1 - 1/M}\frac{(p-1)}{\alpha_pp(p-\alpha p-1)}\left(\sum_{q\in P_M\setminus\{p\}}\frac1{q^2}\right)\ +\  \frac{7 \log_M(n) }{M}.}

By \eqref{alphapboundn}, \eqn{\label{quantityrl1}\sum_{q\in P_0}\lambda_1(q)
&< \frac{ (1 + 1/M)^2}{ 1 - 1/M} \frac{(P(2)-1/4)}{ 0.24} \frac{ (p-1)( 1 -3\alpha/2)}{( p - \alpha p -1 )}
+ \frac{7 \log_M(n) }{M}  \\
&\le \frac{ (1 + 1/M)^2}{ 1 - 1/M} \frac{(P(2)-1/4)}{ 0.24} 
+ \frac{7 }{ \log\log(n) } \\ & < 1} 
for $n$ sufficiently large ($ \log\log(n) \ge 45$ suffices), and we have the desired 2-coprime pair.\\[0mm] 

\noindent
{\bf Case (b)}: For all $p\in P_M$, \eqn{\label{alphapbound}\alpha_p \le \frac{0.24}{p(1-3\alpha/2)}.}

\noindent
For every $p\in P_0$, consider the quantity \eq{\lambda_0(p)=\frac{|S\cap p\Z|}{|S|}\frac{|T\cap p\Z|}{|T|}.} We show that the sum of these terms is less than 1. We bound this term in two cases:
\begin{itemize}
    \item \textbf{$p\in P_M$.} In this case, by \eqref{test} and \eqref{ijpest}, \eqn{\label{l0plm}\lambda_0(p)\le\frac{\alpha_p}\alpha \frac{|J\cap p\Z|}{(1-\alpha)m}<\frac{\alpha_p}\alpha \frac{(1 + 1/M) m/p}{(1-\alpha)m}=\frac{( 1 + 1/M)\alpha_p} {\alpha(1-\alpha)p}.}
    \item \textbf{$p\in P_0\setminus P_M$.} By \eqref{absrem}, (assuming $M\ge15$,) \eqn{\label{tppgm}\frac{|T\cap p\Z|}{|T|}< \frac{m/p+1}{(1-\alpha)m}\le\frac{2}{M}\left(\frac Mp+\frac1{M^2}\right)<\frac{2.01}{M}.}
\end{itemize}
Therefore, by \eqref{l0plm}, \eqref{tppgm}, \eqref{ssps} and \eqref{alphapbound} (noting that $\alpha_p\le \alpha$ and assuming $ M \ge 16$), we have \eqn{\sum_{p\in P_0}\lambda_0(p)&<\sum_{p\in P_M}\frac{17\alpha_p} {16\alpha(1-\alpha)p}\ +\left(\sum_{p\in P_0\setminus P_M}\frac{|S\cap p\Z|}{|S|}\right)\frac{2.01}{M}\nonumber\\
&< \sum_{p\in P_M}\frac{ (1 + 1/M)\alpha_p} {\alpha(1-\alpha)p}\ +\ \log_M(n) \cdot \frac{2.01}{M}\nonumber\\
&\le \label{quantityrl2} (1 + 1/M) \frac{ 0.24}{\alpha( 1 -\alpha)( 1 -3\alpha/2) } \left(\sum_{p\in P_M} \frac{1}{p^2}\right) + \frac{2.01}{\log\log(n)}.}

Define function $f: \alpha \mapsto \alpha( 1 -\alpha)( 1 -3\alpha/2)$, 
and observe that $f$ is concave on the interval $ (0,5/9)$.  It follows that $f$ takes its minimum value in the interval $ [1/16, 1/2]$ at one of the endpoints.  As
$f(1/2) = 1/16$ and  $f(1/16) = 435/2^{13} > 1/19$, 
we have
 \eqn{\sum_{p\in P_0}\lambda_0(p) < (1 + 1/M) 0.24 \cdot 19 (P(2) - 0.25) +  \frac{2.01}{\log\log(n)} < 1,}
for $n$ sufficiently large.  (Here $ \log\log(n) \ge 28$ suffices.)
%
\ep

\vskip5mm
\begin{rem} \label{largen} The explicit conditions on $n$ that are sufficient for the proofs of Lemmas~\ref{rgcoprime}~and~\ref{rlcoprime} play a role when we apply these Lemmas in Section~4. Note that we can establish better bounds by writing some of the conditions in terms of both $n$ and the parameter $M$ (instead of simply applying the bound $ M > \log(n)$). Indeed, the conditions 
\[ M = \left( \frac{m}{5} \right)^{ 1/ \log_2(2r)} > \log n \ \ \ \text{ and } \ \  \ \frac{ \log(n)}{M \log(M)} < \frac{1}{9} \]
are sufficient for Lemma~\ref{rgcoprime}, and the conditions
\[ M = m^{1/3} \ge 16 \ \ \ \text{ and } \ \ \ 
\frac{ \log(n)}{ M \log(M)} < \frac{1}{45} \]
are sufficient for Lemma~\ref{rlcoprime}.
\end{rem}

\vskip5mm


\bp[Proof of Theorem \ref{evencase}]
The four outcomes are clearly pairwise disjoint. We will show that at least one of them happens.

Let $I_1=I\setminus 2\Z$, $I_2=I\cap 2\Z$, $J_1=J\setminus 2\Z$, $J_2=J\cap 2\Z$. They are all integer APs with cardinality $m$ and common difference 2. For $i=1,2$, define $S_i=S\cap I_i$, $T_i=T\cap J_i$. 

Since $|S_1|+|S_2|+|T_1|+|T_2|\ge 2m$, at least one of the following happens:
\begin{itemize}
    \item Case I: $|S_1|+|T_2|>m$ (which implies $S_1,T_2\neq\emptyset$);
    \item Case II: $|S_2|+|T_1|>m$ (which implies $S_2,T_1\neq\emptyset$);
    \item Case III: $|S_1|+|T_2|=|S_2|+|T_1|=m$.
\end{itemize}

If $|S_1|+|T_2|\ge m$ and $S_1,T_2\neq \emptyset$, then by Lemma \ref{interlem}, there are $s\in S_1$, $t\in T_2$ that are 2-coprime. Because $s$ is odd, that means they are coprime, so the last outcome applies. Analogously, a coprime pair exists if $|S_2|+|T_1|\ge m$ and $S_2,T_1\neq\emptyset$. To \textit{not} fall into the last outcome, we must be in Case III, with at least one of $S_1,T_2$ empty and at least one of $S_2,T_1$ empty.

\begin{itemize}
    \item If $S_1,S_2=\emptyset$, we are in Outcome 1.
    \item If $T_2,T_1=\emptyset$, we are in Outcome 2.
    \item If $S_1,T_1=\emptyset$, we are in Outcome 3.
    \item If $T_2,S_2=\emptyset$, then $|S_1|=|T_1|=m$. By Lemma \ref{interlem} there are $s\in S_1$, $t\in T_1$ that are 2-coprime. Because both $s$ and $t$ are odd, they are coprime, so we are in Outcome 4.
\end{itemize}
In conclusion, one of the four outcomes must happen, so we are done.
\ep

\section{The lonely runner with slow runners}
In this Section
we prove Theorem \ref{thm:run}.

Let $\C_0$ be the constant in Theorem \ref{evencase}. Set $\C$ as a constant for which the following holds for sufficiently large $n$: \eq{\exp(\C\cdot(\log\log n)^2)\ge 8\ceil{\exp(\C_0\cdot(\log\log(2n))^2)}+6.} For convenience, let \eq{k=k(n)=\frac{\exp(\C\cdot(\log\log n)^2)}2,\ \M=\M(n)=\ceil{\exp(\C_0\cdot(\log\log(2n))^2)}.}
The assumption then becomes $k\ge 4\M+3$. We will show that for sufficiently large $n$, if positive integers $v_1<v_2<\dots<v_n$ satisfy that $n<v_n\le 2n-2k$, then Condition \eqref{loose} holds. 

Consider such $v_1,\dots,v_n$. Denote $V=\{v_1,\dots,v_n\}$. Since $v_n>n$, in particular $V\neq[n]$, so there exists a ``largest missing number in $[n]$'', \eq{x=\max\left([n]\setminus V\right).} 

We first claim that if $x>n-k$ then Condition \eqref{loose} holds. Note that $x\neq v_i$ for all $i$, and larger multiples of $x$ are $>2n-2k$, which is too large to be one of $v_1,\dots,v_n$. Letting $t=1/x$, the quantity $\{v_it\}=\{v_i/x\}$ cannot be zero; moreover, its denominator is at most $x\le n$. Thus, $1/(n+1)<\{v_it\}<n/(n+1)$ for all~$i$.

Hence, we may assume $x\le n-k$. Next, we claim that it is possible to cut $[n]$ into some groups of consecutive integers in ascending order, such that all but the last group has cardinality either $2\M$ or $2\M+2$, the last group's cardinality is between $2\M+2$ and $4\M+3$, and $x$, $x+1$ and $x+2$ belong to the same group. One can start with $\{1,\dots,2\M+2\},\{2\M+3,\dots,4\M+4\}$, and continue to append groups of size $2\M+2$ by default; but if a new group would cut $x$, $x+1$ and $x+2$ apart, shrink its size by 2 to avoid the separation. There will likely be a residue of size $<2\M+2$ at the end of $[n]$, in which case merge it with the previous group. Given that the distance between $x$ and $n$ is at least $k\ge 4\M+3$, the resulting last group would not contain $x$, hence unaffected by the potential shrinking.

Denote the groups $I_1,I_2,\dots,I_\ell,I_{\ell+1}$ in ascending order, and say $x,x+1,x+2\in I_r$ for some $r\in[\ell]$. Now we will also partition $\{n+1,\dots,2n\}$ into groups. Let $J_\ell$ be the interval starting at $n+1$ with the same cardinality as $I_\ell$. Let $J_{\ell-1}$ be the next contiguous interval, with the same cardinality as $I_{\ell-1}$. (That is, $J_{\ell-1}=\{n+|I_\ell|+1,n+|I_\ell|+2,\dots,n+|I_\ell|+|I_{\ell-1}|\}$.) Define $J_{\ell-2},\dots,J_1$ analogously. The remaining numbers would form the interval $J_{\ell+1}$, whose cardinality is equal to that of $I_{\ell+1}$.
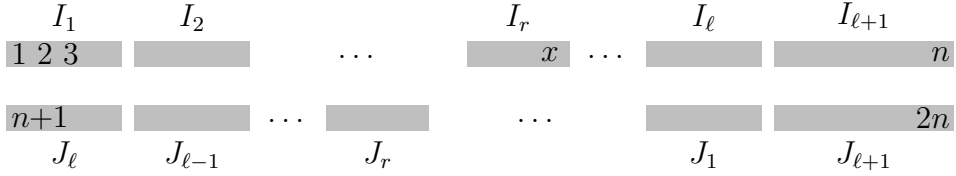
\begin{center}
\begin{tikzpicture}[scale=0.17]
\fill [lightgray] (0,0) -- (9,0) -- (9,-2) -- (0,-2);
\draw (4.5,0) [above] node {$I_1$};
\fill [lightgray] (10,0) -- (19,0) -- (19,-2) -- (10,-2);
\draw (14.5,0) [above] node {$I_2$};
\fill [lightgray] (36,0) -- (44,0) -- (44,-2) -- (36,-2);
\draw (40,0) [above] node {$I_r$};
\fill [lightgray] (50,0) -- (59,0) -- (59,-2) -- (50,-2);
\draw (54.5,0) [above] node {$I_\ell$};
\fill [lightgray] (60,0) -- (74,0) -- (74,-2) -- (60,-2);
\draw (67,0) [above] node {$I_{\ell+1}$};
\draw (-0.5,-1) [right] node {$1\ 2\ 3$};
\draw (27.5,-0.2) [below] node {$\dots$};
\draw (44,-1) [left] node {$x$};
\draw (47,-0.2) [below] node {$\dots$};
\draw (74.6,-1) [left] node {$n$};

\fill [lightgray] (0,-5) -- (9,-5) -- (9,-7) -- (0,-7);
\draw (4.5,-7) [below] node {$J_\ell$};
\fill [lightgray] (10,-5) -- (19,-5) -- (19,-7) -- (10,-7);
\draw (14.5,-7) [below] node {$J_{\ell-1}$};
\fill [lightgray] (25,-5) -- (33,-5) -- (33,-7) -- (25,-7);
\draw (29,-7) [below] node {$J_r$};
\fill [lightgray] (50,-5) -- (59,-5) -- (59,-7) -- (50,-7);
\draw (54.5,-7) [below] node {$J_1$};
\fill [lightgray] (60,-5) -- (74,-5) -- (74,-7) -- (60,-7);
\draw (67,-7) [below] node {$J_{\ell+1}$};
\draw (-0.5,-6) [right] node {$n\hspace{-0.04in}+\hspace{-0.04in}1$};
\draw (22,-5.2) [below] node {$\dots$};
\draw (41.5,-5.2) [below] node {$\dots$};
\draw (74.6,-6) [left] node {$2n$};
\end{tikzpicture}
\end{center}

Note inductively that for all $1\le j\le \ell$, $\min (I_j)+\max(J_j)=2n+1-|I_{\ell+1}|$. Thus, \eq{\min(I_j)+\min(J_j)&=2n+1-|I_{\ell+1}|-(|J_j|-1)\\
&\ge 2n+1-(4\M+3)-(2\M+1)\\
&=2n-6\M-3\\
&>2n-2k.\\
\max(I_j)+\max(J_j)&=2n+1-|I_{\ell+1}|+(|J_j|-1)\\
&\le 2n+1-(2\M+2)+(2\M+1)\\
&=2n.}

For all $1\le j\le \ell+1$, let $S_j=I_j\setminus V$ and $T_j=J_j\setminus V$. We claim that if for some $j\in[\ell]$ there exist $s\in S_j$ and $t\in T_j$ coprime then Condition \eqref{loose} holds. By the above calculation, $2n-2k< s+t\le 2n$. Consider $t=q/(s+t)$, where $q$ is the inverse of $s$ in $(\Z/(s+t)\Z)^\times$. Observe that no $v_i$ can be $0$, $s$ or $t$ modulo $s+t$, because $s,t$ are non-speeds by definition and $s+t>2n-2k$ is too large. But these are the only ways to have $v_iq$ be $0,1$ or $-1$ modulo $s+t$. Thus, for all $i$, $v_iq/(s+t)$ is at least $2/(s+t)\ge 2/(2n)$ away from the nearest integer. Hence $1/(n+1)<\{v_it\}<n/(n+1)$.

Define $\alpha_j=|S_j|+|T_j|$ and $m_j=|I_j|/2=|J_j|/2$ for all $1\le j\le \ell+1$. By the grouping conditions, for all $j\le \ell$ we have $m_j\in \Z$ and $m_j\ge \M$. Note that $\sum_{j=1}^{\ell+1}\alpha_j=n=\sum_{j=1}^{\ell+1}2m_j$. Recall that $I_{\ell+1}\subseteq V$ because $x\le n-k\le n-|I_{\ell+1}|$, and $J_{\ell+1}\cap V=\emptyset$ because $v_n\le 2n-2k\le2n-|J_{\ell+1}|$. Hence, $\alpha_{\ell+1}=0+|J_{\ell+1}|=2m_{\ell+1}$. By pigeonhole, either $\alpha_j>2m_j$ for some $j\in[\ell]$, or $\alpha_j=2m_j$ for all $j\in[\ell]$. In the former case, according to Theorem~\ref{evencase} applied on $(2n,m_j,I_j,J_j,S_j,T_j)$, there must be a coprime pair between $S_j$ and $T_j$, witnessing Condition \eqref{loose}. In the second case, consider $S_r$ and $T_r$. We have $|S_r|+|T_r|=\alpha_r=2m_r$, so by Theorem \ref{evencase}, one of the four outcomes happens: $S_r=\emptyset$, $T_r=\emptyset$, ($S_r=I_r\cap 2\Z$ and $T_r=J_r\cap 2\Z$), or there is a coprime pair between $S_r$ and $T_r$. But as $x\notin V$, $S_r\neq\emptyset$; as $x+1\in V$, $T_r\neq\emptyset$; as $x+1$ and $x+2$ have different parities and are both in $V$, $S_r$ misses a number of each parity. Thus, only the last outcome is possible, which also witnesses Condition \eqref{loose}.

\section{Coprime mappings for adjacent intervals}

In this Section we prove Theorem~\ref{tnconj}; that is, we show that if $ 0 \le \ell < k$ and $ k \ge 4 $ then there is coprime mapping from $A=\{\ell+1,\dots,\ell+k\}$ to $\ B=\{\ell+k+1,\dots,\ell+2k\}$.

\noindent

Theorem~\ref{tnconj} follows from Lemmas~\ref{rgcoprime}~and~\ref{rlcoprime} when $n$ is sufficiently large. 
We prove Theorem~\ref{tnconj} for smaller values of $n$ using a separate argument. We note in passing that 
there are parallels between this argument and the earlier 
work of Pomerance and Selfridge \cite{pomerance1980proof}.  Both arguments make use of a statistical study of the parameter $ \phi(x)/x$, where $\phi$ is the Euler totient function, and both proofs make use of estimates on the prime counting functions given by Rosser and Schoenfeld \cite{older} (but our reliance on these estimates is significantly less extensive).

To handle smaller $n$ we need some additional definitions and elementary
observations.  Let $ p_i$ be the $i^{\rm th}$ odd prime, and set $q_i = \prod_{j=1}^i p_i$. For each integer $x$ we let $ P(x) $ be the set of odd prime factors of $x$ and define
\eq{ \gamma(x) = \prod_{p \in P(x)} \frac{ p-1}{p}.}
Note that $ \gamma(x)$ is the proportion of numbers that are 2-coprime with $x$.  (Note further that $ \gamma(x) = \phi(x)/x$ when $x$ is odd.)
\begin{cla}\label{prodarg}
Let $J$ be an integer AP with cardinality~$m$ and common difference~$d$. Let $x$ be an integer such that $d$ is coprime with all elements of $P(x)$. Then more than $\gamma(x) m - 2^{|P|}+1$ numbers in $J$ are 2-coprime with $x$.
\end{cla}
\bp
For $p\in P$, let $J_p=J\cap p\Z$. By inclusion-exclusion, \eq{\left|\bigcup_{p\in P} J_p\right|=\sum_{\emptyset\neq Q\subseteq P} (-1)^{|Q|+1}\left|\bigcap_{q\in Q}J_q\right|.}

Because $P\cup\{d\}$ is mutually coprime, the cardinality of $\bigcap_{q\in Q}J_q$ can be approximated: \eq{\left|\bigcap_{q\in Q}J_q\right|-\frac{m}{\prod_{q \in Q}q}\in (-1,1).}
Note that $m/\prod_{q \in Q}q$ is the ``heuristic'' cardinality as the size of $J$ gets large. 
By summing these terms for all nonempty subsets of $P$, the cardinality of $\bigcup_{p\in P} J_p$ can be approximated, \eq{\left(1-\prod_{p\in P} \frac{p-1}p\right)m=(1-\gamma) m,} with an error less than $2^{|P|}-1$. The claim follows by taking the complement.
\ep
\begin{cla}
\label{single}
Let $J \subset [n]$ be an AP with common difference 1 or 2 such that $ |J| =m \ge 2$. Let $ T \subseteq J$ such that $ |T| \ge (m+1)/2$. If $ x $ is an integer such that $ |P(x)| \le 1$ then there is $ y \in T$ such that $x$ and $y$ are 2-coprime.
\end{cla}
\bp
Since $ m \ge 2$, the set $T$ contains two consecutive elements or two elements out of 3 consecutive elements of $J$. In either case, at least one of these numbers is not divisible by the odd prime that divides $x$.
\ep


With these preliminary observations in hand we are ready to prove Theorem~\ref{tnconj}.
  Set $ A_0 = A \cap 2 \Z$, $A_1 = A \setminus A_0$,  $ B_0 = B \cap 2 \Z$, $B_1 = B \setminus B_0$.  Note that it suffices to find a 2-coprime mapping from $ A_0$ to $B_1$ and another from $A_1$ to $B_0$.  (As the intervals $A$ and $B$ are consecutive we have $|A_0| = |B_1|$ and $|A_1| = |B_0|$.)  Hall's Theorem states that there is a 2-coprime mapping between  $ A_1$ and $B_0$ if and only if for every pair of sets $ S \subseteq A_1$ and $ T \subseteq B_0 $ such that $ |S| + |T| = |A_0| + 1 = |B_1| + 1$ there exists a 2-coprime pair $x,y$ such that $x\in S$ and $y\in T$. When $k$ is sufficiently large we can apply Lemmas~\ref{rgcoprime}~and~\ref{rlcoprime} -- with $ n = 2k+\ell$ and $m = |A_0| = |B_1|$ or $ m = |A_1|=|B_0|$ -- to conclude that the desired 2-coprime pair exists, and hence the desired 2-coprime mappings exist.

Let  $I$ and $J$ be disjoint APs in $[n]$ with common difference 2 such
that $|I| = |J| = m \ge \frac{n-2}{6}$. (We take $ \{ I,J \} = \{ A_0, B_1\} $ or $ \{I,J\} = \{ A_1, B_0\}$. Note that $ \frac{ n-2}{6} \le \lfloor k/2 \rfloor $.)
Let $ S \subseteq I$ and $T \subseteq J$ such that $ |S| + |T| = m+1$ and $ |S| \le |T|$.  As above we set $r = |S|/m$.  We show that there is a 2-coprime pair $x, y$ such that $ x \in S$ and $y \in T$. We consider two cases depending on the value of $n$.  For large $n$ we appeal to Lemmas~\ref{rgcoprime}~and~\ref{rlcoprime}, and for small $n$ we provide a direct argument.\\[-1mm]

\noindent
{\bf Case 1:} $ 2k+ \ell = n > 3 \cdot 10^7$.\\[-1mm]

First consider $ 2 \le r \le 16$.  Here we apply Lemma~\ref{rlcoprime}.  We clearly 
have $ m > \log(n)^3$.  As noted in Remark~\ref{largen}, the parameter $n$ is sufficiently large if the following two conditions hold:
\[ m^{1/3} = M \ge 16\ \ \ \text{ and } \ \ \ \frac{ \log(n)}{\frac{1}{3} \log(m) \cdot m^{1/3} } = \frac{ \log_M(n)}{ M} < \frac{1}{45}. \]
Both conditions hold in the range in question, 
and Lemma~\ref{rlcoprime} gives the desired 2-coprime pair.

Now consider $ r > 16$. We first consider $ 3 \cdot 10^7 < 
n < 10^{50}$.
Assume for the sake of contradiction that we do not have the desired 2-coprime pair. Then it follows from Claim~\ref{prodarg} that for any $ x \in S$
we have
\eq{ 15m/16 < m -|S| < |T| < m - m \gamma(x) + 2^{|P(x)|}.}
Observe that $ |P(x)| \le \log_3(n) $ for all integers $x$.  Furthermore, as $ q_{32} > 10^{50}$, the assumed restriction on $n$ implies $ \gamma(x) > \gamma( q_{31} ) > 1/5$.  
It follows that we have
\begin{multline*}
1/5 < \gamma(x) < 1 /16 + 2^{ |P(x)|}/m < 1/16 + 2^{ \log_3(n)}/( (n-2)/6)  \\
< 1/16 + 6 \cdot n^{ \frac{\log(2)}{ \log(3)} -1} + 12/n < 1/5,  
\end{multline*}
which is a contradiction.

It remains to consider $ r > 16$ and $ n > 10^{50}$.  
We apply Lemma~\ref{rgcoprime} when $r < \log(n)^2$. In order to handle the
large $r$ case, we first note that if $ r > \log(n)^2$ then $T$ either contains 
primes $ p_1, p_2 > k/4$ or contains $ 2 p_1, 2p_2$ where $ p_1, p_2 > k/8$ are primes.  To see this, we appeal to bounds on the prime counting function $ \pi(x)$ (see Theorem~1 in \cite{older}).  First suppose $J$ consists of odd numbers and $a$ is the largest element of $J$. Then the number of primes in $ J \cap [a-m+1,a]$ is
\begin{equation*}
\pi(a) - \pi( a - m ) \ge \frac{a}{ \log(a)} 
-  \frac{a - m }{ \log(a - m)}  \left(1+\frac2{\log (a - m)}\right)
\end{equation*}
For ease of notation, let $ a - m = \eta a$ we have
\eq{
\pi(a) - \pi( \eta a ) & \ge \frac{a}{ \log(a)} 
-  \frac{\eta a }{ \log(\eta a)}  \left(1+\frac2{\log (\eta a)}\right) \\
& = \frac{ a - \eta a}{ \log(a)} + \frac{ \eta \log( \eta) \cdot a}{ \log( \eta a) \log (a) } 
- \frac{ 2 \eta a}{ \log( \eta a)^2} \\
& \ge \frac{m}{ \log(n)} - \frac{2n}{ \log(n)^2} \\
& \ge \frac{2 m}{ \log(n)^2 } \\
& > |S|.
}
Recalling that $ |S| + |T| \ge m + 1$, the number of elements of $J$ that are NOT elements of $T$ is at most $ |S|-1$.  Thus we have the two desired primes in $T$.   If $J$ consists of even numbers then we apply the same estimates to $ J/2 = \{ x/2: x \in J \}$ to get the desired elememts $ 2p_1, 2p_2$ where $ p_1$ and $p_2$ are prime. In either case,
no element of $S$ is divisible by both $ p_1$ and $ p_2$, and we have the desired 2-coprime pair.

When $16 < r < \log(n)^2 $ and $ n > 10^{50}$ we apply Lemma~\ref{rgcoprime}.  Recall that, as noted in Remark~\ref{largen}, $n$ is sufficiently large if the following conditions hold:
\[ M = \left( \frac{m}{5} \right)^{ 1/ \log_2(2r)} > \log n \ \ \ \text{ and } \ \  \ \frac{ \log(n)}{M \log(M)} < \frac{1}{9}. \]
These conditions hold here (indeed, we have $ M > 3 \log(n)/2$), and Lemma~\ref{rgcoprime} applies.\\[-1mm]

\noindent{\bf Case 2:} $ 83 \le 2k + \ell = n < 3 \cdot 10^7$ \\[-1mm]

\noindent
Here we apply Claim~\ref{prodarg}. Note if
 $ x \in S$ and
 \eqn{\label{gmajor}  |T| \ge ( 1 -\gamma(x)) m + 2 ^{|P(x)|} - 1,}
then Claim~\ref{prodarg} implies that the desired 2-coprime pair exists. As 
many elements of $ [ 3 \cdot 10^7]$ have large values of $ \gamma(x) $, this 
observation is usually sufficient to complete the proof for $n$ in this 
interval. In order to make the proof
precise we consider cases. In some cases we will need to make a more detailed study of the collection of sets $ \{ P(x): x \in S\}$.\\[-1mm]

\noindent
{\bf Case 2a:} There exists $s\in S$ that is not divisible by 3. \\[-1mm]
 
 \noindent
As $ |T| \ge (m+1)/2 $ it suffices to show that there exists $ x \in S$ such that
\eqn{ \label{no3}  (\gamma(x) - 1/2) m \ge 2 ^{|P(x)|} - 3/2  }
Set $ w_a = q_{a+1}/3 = \prod_{i=2}^{a+1} p_i$.
Among numbers $x$ such that $ 3 \nmid x$ with a fixed value of $a= |P(x)|$, $x$ is minimized by $w_a$ and the parameter $ \gamma(x)$ is minimized by $ \gamma( w_a )$. Set \eq{\chi_a=\frac{2 ^a - 3/2}{\gamma(w_a) - 1/2}\cdot 6+2} and observe that if for a particular value of $a$ we have $n\ge \chi_a$ 
and there exists $ s \in S\setminus 3\Z$ such that $ |P(s)| = a $
then\footnote{We implicitly also need $\gamma(w_a)>1/2$, which holds within our range of consideration.} we have (\ref{no3}) and Claim~\ref{prodarg} implies that the desired 2-coprime pair exists.
We refer to the following table for values of $w_a$ and $\chi_a$. 
\begin{center}
\begin{tabular}{|r|r|r|r|}
\hline
$a$ & $w_a$  & $\gamma( w_a)$ & $\chi_a$ \\
\hline 1 & 5 & 0.8 & 12 \\
\hline 2 & 35 & 0.6857 & 82.8 \\
\hline 3 & 385 & 0.6234 & 318.1 \\
\hline 4 & 5,005 & 0.5754 & 1155.5 \\
\hline 5 & 85,085 & 0.5416 & 4403.6 \\
\hline 6 & 1,616,615 & 0.5131 & 28689.1 \\
\hline
\end{tabular}
\end{center}
Take any $s\in S\setminus 3\Z$ and let $a=|P(s)|$. As $w_7>3\cdot 10^7$, we have $a\le 6$. If $a\le 2$, we have $n\ge 83>\chi_a$ immediately. If $a\ge 3$, then we also have $n\ge s\ge w_a \ge \chi_a$. Thus, we always have the desired 2-coprime pair.\\[-1mm]

\noindent
{\bf Case 2b:} Every number in $S$ is divisible by 3. \\

\noindent
As $|S|\le |I\cap 3\Z|\le(m+2)/3$ and hence $|T|\ge(2m+1)/3$, 
it suffices to show that there exists $ x \in S$ such that
\eqn{ \label{with3}  (\gamma(x) - 1/3) m \ge 2 ^{|P(x)|} - 4/3.  }
Proceeding as in the previous case, if there exists a value of $a$ such that
\eq{n\ge \kappa_a := \frac{2 ^a - 4/3}{\gamma(q_a) - 1/3}\cdot 6+2}
and $ x\in S$ such that $ |P(x)| =a$ then Claim~\ref{prodarg} implies that the desired
2-coprime pair exists.  It follows from the following table that we have the desired condition
for $ a=1,2,4,5,6,7$. 
\begin{center}
\begin{tabular}{|r|r|r|r|}
\hline
$a$ & $q_a$  & $\gamma( q_a)$ & $\kappa_a$ \\
\hline 1 & 3 & 0.6667 & 14 \\
\hline 2 & 15 & 0.5333 & 82 \\
\hline 3 & 105 & 0.4571 & 325.1 \\
\hline 4 & 1,155 & 0.4156 & 1071.9 \\
\hline 5 & 15,015 & 0.3836 & 3661.3 \\
\hline 6 & 255,255 & 0.3611 & 13567.5 \\
\hline 7 & 4,849,845 & 0.3420 & 87210.9 \\
\hline
\end{tabular}
\end{center}
As $ q_8 > 3 \cdot 10^7$ we do not need to consider larger values of $a$. It remains to consider the case where $ |P(x)| = 3$ for all $ x \in S$. Again, if $n \ge \kappa_3=325.1$ we will have the desired 2-coprime pair, so assume $n\le 325$. 
Note that \eq{S\subseteq \{x\in [325]:|P(x)|=3\}=\{105,165,195,210,231,255,273,285,315\}.} For \eqref{gmajor} to not hold, we must have \eq{m + 1 - |S| \le |T| &< (1-\gamma(q_3))m+2^3-1\\\Rightarrow m&<\frac{6+|S|}{\gamma(q_3)}\\\Rightarrow \max(S)\le n&< 6\cdot\frac{6+|S|}{\gamma(q_3)}+2<14|S|+81.} One can verify that this is not possible. Hence, \eqref{gmajor} holds and we have the desired 2-coprime pair.\\[-1mm]

\noindent{\bf Case 3:} $ 8 \le 2k + \ell = n < 83$ \\[-1mm]

\noindent
Here we make a more careful analysis of the collections of sets $ \{ P(x): x \in S \} $. We recall Claim~\ref{single}: If any of these sets has cardinality at most 1 then we have the desired 2-coprime pair. It follows that, as $ n<83< q_3$, we may assume that $ |P(x)|=2$ for all $ x \in S$, and so we can view
 $ \{ P(x): x \in S \} $ as a graph $G$ on a vertex set consisting of the odd primes.  
Now we consider some further cases.\\

\noindent{\bf Case 3a:}  $G$ contains disjoint edges $p_1p_2$ and $ p_3p_4$.

\begin{quote}
Note that in this case we have $n \ge 35$ and $ m \ge 6$. Note that a number $y$ is 2-coprime with neither of the corresponding elements of $S$ if and only if $ P(y) $ contains a set in $ \{ p_1, p_2\} \times \{ p_3, p_4\}$.
As 15 is the smallest such product, at most 4 elements out of any 15 consecutive elements of $J$ contains one of these products. 
It follows that we have the desired 2-coprime pair so long as $m \ge 8$. In the case $ 4 \le m \le 7$ we have $ n \le 44$ and there are at most 4 elements $z$ of $ I \cup J$ such that $ |P(z)| \ge 2$. Two of them are already in $S$, and 
$T$ has at least 3 elements, so some element $t\in T$ satisfies $|P(t)|\le 1$, giving the desired 2-coprime pair.  We finally note that if $ m \le 3$ then $ n \le 20$ and this case is not possible.    
\end{quote}

%

\noindent{\bf Case 3b:} $G$ has a vertex of degree 2.

\begin{quote}
    Suppose $G$ contains the edges $ p_1p_2$ and $p_1p_3$.  Note if $y$ is an element of $J$ such that $ p_1 \not\in P(y)$ and $ \{p_2,p_3\} \not\subseteq P(y)$ then $y$ is 2-coprime with one of the corresponding elements of $S$.  
    
    First consider $ p_1 > 3$.  In this case $ n \ge 35$, which implies $ m \ge 6$. Furthermore, if $ p_1> 3$ then at most 3 out of 10 consecutive elements of $J$ are not 2-coprime with the corresponding elements of $S$.  It follows from these observations that we have the desired 2-coprime pair.
    
    So, we may assume $ p_1=3$. As $m\le (k+1)/2 \le (n+2)/4 < 22$, $J$ contains at most one element that is a multiple of $ p_2p_3$, and the number of elements in $J$ that are not 2-coprime to any of the corresponding elements of $S$ is at most $ (m+2)/3 + 1$. As $|T| \ge (m+1)/2$, we have the desired 2-coprime pair if $ m \ge 8$.  
    
    Finally consider $m \le 7$, which implies $ n \le 44$.  Here we claim that $J$ does not contain any element that is a multiple of $ p_2p_3$.  Assume for the sake of contradiction that $J$ has such an element. Given the bound on $n$, 35 is the only product of distinct odd primes in $ I \cup J$ that does not include 3. So, we have $ 35 \in J$. Recalling that one of $I$ and $J$ consists of only even numbers, we observe that $ 35 \in J$ implies that $I$ consists of even numbers and $ 30, 42 \in S$. This is a contradiction as $A$ and $B$ are disjoint intervals.  It follows that the number of elements of $J$ that are not 2-coprime to the corresponding elements of $S$ is at most $ (m+2)/3$, and we have the desired 2-coprime pair.
\end{quote}

\noindent{\bf Case 3c:} $G$ consists of a single edge $ p_1p_2$

\begin{quote}
    Here we observe that $ |S| \le 2$ and that at most three elements out of any 5 consecutive elements of $J$ is not 2-coprime with the corresponding elements of $S$.  As $ |T| \ge m +1  -|S| \ge m-1$, we have the desired 2-coprime pair when $ m\ge 5$.  But $m \le 4$ implies $ n \le 26 $, which implies $|S| \le 1$ and $ T =J$. As at least 2 out of any 3 consecutive elements of $J$ are 2-coprime with the corresponding elements of $S$, we have the desired 2-coprime pair when $ m = 3,4$. Finally, if $m\le2$ then we have $ n \le 14$ and this case is not possible. 
\end{quote}

\bibliographystyle{plain}
\bibliography{bib.bib}
\end{document}